\documentclass[12pt]{amsart}
\usepackage{amscd,amsmath,amssymb,amsfonts}
\usepackage[all]{xy}
\theoremstyle{plain}
\newtheorem{thm}{Theorem}
\newtheorem{lem}[thm]{Lemma}
\newtheorem{cor}[thm]{Corollary}
\newtheorem{prop}[thm]{Proposition}

\theoremstyle{definition}
\newtheorem{defn}[thm]{Definition}

\numberwithin{thm}{section} \numberwithin{equation}{section}

\newcommand{\ga}[2]{\begin{gather}\label{#1}#2 \end{gather}}

\newcommand{\surj}{\twoheadrightarrow}
\newcommand{\inj}{\hookrightarrow}

\newcommand{\Spec}{{\rm Spec \,}}

\newcommand{\sE}{{\mathcal E}}
\newcommand{\sF}{{\mathcal F}}

\newcommand{\sO}{{\mathcal O}}

\newcommand{\sV}{{\mathcal V}}

\newcommand{\Q}{{\mathbb Q}}

\newcommand{\Z}{{\mathbb Z}}

\begin{document}

\title{Direct images of bundles under Frobenius morphism}
\author{Xiaotao Sun}
\address{Academy of Mathematics and Systems Science, Chinese Academy of Science, Beijing, P. R. of China}
\email{xsun@math.ac.cn}
\address{}
\date{February 22, 2008}
\thanks{Partially supported by a matched grant of Chinese Academy of Science and
the DFG Leibniz Preis of Esnault-Viehweg}
\begin{abstract} Let $X$ be a smooth projective variety of dimension $n$ over an
algebraically closed field $k$ with ${\rm char}(k)=p>0$ and $F:X\to
X_1$ be the relative Frobenius morphism. For any vector bundle $W$
on $X$, we prove that instability of $F_*W$ is bounded by
instability of $W\otimes{\rm T}^{\ell}(\Omega^1_X)$ ($0\le \ell\le
n(p-1)$)(Corollary \ref{cor3.8}). When $X$ is a smooth projective
curve of genus $g\ge 2$, it implies $F_*W$ being stable whenever $W$
is stable.
\end{abstract}
\maketitle
\begin{quote}
Dedicated to Professor Zhexian Wan on the occasion of his 80th
birthday.
\end{quote}
\section{Introduction}

Let $X$ be a smooth projective variety of dimension $n$ over an
algebraically closed field $k$ with ${\rm char}(k)=p>0$. Fix an
ample divisor ${\rm H}$ on $X$, by a semistable (resp. stable)
torsion free sheaf, we mean a ${\rm H}$-slope semistable (resp.
${\rm H}$-slope stable) sheaf in this paper. For a torsion free
sheaf $\sF$ on $X$, there is a unique filtration
$$0=\sF_0\subset \sF_1\subset\cdots\subset \sF_k=\sF$$
such that $\sF_i/\sF_{i-1}$ ($1\le i\le k$) are semistable torsion
free sheaves and
$$\mu_{{\rm max}}(\sF):=\mu(\sF_1)>\mu(\sF_2/\sF_1)>\cdots>\mu(\sF_k/\sF_{k-1}):=\mu_{{\rm
min}}(\sF).$$ The instability of $\sF$ was defined as ${\rm
I}(\sF)=\mu_{{\rm max}}(\sF)-\mu_{{\rm min}}(\sF)$, which measures
how far from $\sF$ being semi-stable. In particular, $\sF$ is
semi-stable if and only if ${\rm I}(\sF)=0$. On the other hand,
there are sub-bundles ${\rm
T}^{\ell}(\Omega^1_X)\subset(\Omega^1_X)^{\otimes\ell}$,
$0\le\ell\le n(p-1)$, which are the associated bundles of
$\Omega^1_X$ through some elementary (perhaps interesting)
representations of ${\rm GL}(n)$. These representations do not
appear in characteristic zero.

Let $F:X\to X_1$ be the relative Frobenius morphism, for any
vector bundle $W$ on $X$, let ${\rm I}(W,X)$ be the maximal value
of ${\rm I}(W\otimes{\rm T}^{\ell}(\Omega^1_X))$ where $0\le
\ell\le n(p-1)$. Then one of our results in this paper shows
(Corollary \ref{cor3.8}): When $K_X\cdot{\rm H}^{n-1}\ge 0$, we
have
$${\rm I}(F_*W)\le p^{n-1}{\rm rk}(W)\,{\rm I}(W,X)\,.$$
In particular, if the bundles $W\otimes{\rm
T}^{\ell}(\Omega^1_X)$, $0\le \ell\le n(p-1)$, are semistable,
then $F_*W$ is semistable. In fact, when $K_X\cdot{\rm
H}^{n-1}>0$, we can show that the stability of $W\otimes{\rm
T}^{\ell}(\Omega^1_X)$, $0\le \ell\le n(p-1)$, implies the
stability of $F_*W$ (Theorem \ref{thm3.7}).

The main theorem has an immediate corollary that when $X$ is a
smooth projective curve of genus $g\ge 2$, the stability of $W$
implies stability of $F_*W$. This is in fact our original motivation
stimulated by a question raised by Herbert Lange at a conference.
When $W$ is a line bundle, it is due to Lange and Pauly
(\cite[Proposition 1.2 ]{LP}). The present version is based on our
earlier preprint (\cite{Su}), where the theorem was completely
proved only for curves. It should be pointed out, in case of curves,
Mehta and Pauly have proved independently that semi-stability of $W$
implies semi-stability of $F_*W$ by a different method. However,
their method was not able to prove the stability of $F_*W$ when $W$
is stable. In fact, they asked the question: Is stability also
preserved by $F_*$ ? (cf. \cite[Section 7]{MP} for the discussions).

To describe the idea of proof, let us compare it to its opposite
case, a Galois {\'e}tale $G$-cover $f:Y\to X$. Recall that for a
semi-stable bundle $W$ on $Y$, to prove semistability of $f_*W$,
one uses the fact that $f^*(f_*W)$ decomposes into pieces of
$W^{\sigma}$ ($\sigma\in G$). To imitate this idea for $F:X\to
X_1$, we need a similar decomposition of $V=F^*(F_*W)$. Indeed,
use the canonical connection $\nabla: V\to V\otimes\Omega^1_X$,
Joshi-Ramanan-Xia-Yu have defined in \cite{JRXY} for ${\rm
dim}(X)=1$ a canonical filtration
$$0=V_p\subset V_{p-1}\subset\cdots\subset V_{\ell}\subset
V_{\ell-1}\subset\cdots V_1\subset V_0=V$$ such that
$V_{\ell}/V_{\ell+1}\cong W\otimes(\Omega^1_X)^{\otimes \ell}$. It
is this filtration and its generalization that we are going to use
for the study of $F_*W$.

As the first step, we generalize the canonical filtration to higher
dimensional $X$. Its definition can be generalized straightforwardly
by using the canonical connection $\nabla:V\to V\otimes\Omega^1_X$.
The study of its graded quotients are much involved. We show
(Theorem \ref{thm2.7}) that there exists a canonical filtration
$$0=V_{n(p-1)+1}\subset V_{n(p-1)}\subset\cdots\subset V_1\subset
V_0=V=F^*(F_*W)$$ such that $\nabla$ induces injective morphisms
$V_{\ell}/V_{\ell+1}\xrightarrow{\nabla}
(V_{\ell-1}/V_{\ell})\otimes\Omega^1_X$ of vector bundles and the
isomorphisms $V_{\ell}/V_{\ell+1}\cong W\otimes {\rm
T}^{\ell}(\Omega^1_X)$, where ${\rm T}^{\ell}(\Omega^1_X)\subset
(\Omega^1_X)^{\otimes \ell}$ are subbundles given by representations
of ${\rm GL}(n)$ (cf. Definition \ref{defn2.4}). In characteristic
zero, ${\rm T}^{\ell}(\Omega^1_X)={\rm Sym}^{\ell}(\Omega^1_X)$. In
characteristic $p>0$, ${\rm T}^{\ell}(\Omega^1_X)\cong {\rm
Sym}^{\ell}(\Omega^1_X)$ only for $\ell<p$. In general, there is a
resolution of ${\rm T}^{\ell}(\Omega^1_X)$ (Proposition
\ref{prop2.5}) by symmetric powers of $\Omega^1_X$ and exterior
powers of $F^*(\Omega^1_X)$ (After \cite{Su} appeared, Indranil
Biswas told me that a similar filtration was defined and studied in
Proposition 4.1 of their preprint \cite{BH}. However, since their
map (4.7) was wrong, the Proposition 4.1 (also Proposition 4.2
consequently) of \cite{BH} was wrong. After we pointed out these
gaps, they have corrected these mistakes in \cite{BH2}).

To prove the main theorem, we also need to compare sub-sheaves of
$V_{\ell}/V_{\ell+1}$ to sub-sheaves of
$V_{n(p-1)-\ell}/V_{n(p-1)-\ell+1}$ which are $\nabla$-invariant
(Proposition \ref{prop3.6}). It is reduced to consider the (graded)
$K$-algebra
$$R=\frac{K[y_1,y_2,\cdots,y_n]}{(y_1^p,y_2^p,\ldots,y_n^p)}=
\bigoplus^{n(p-1)}_{\ell=0}R^{\ell}$$ with a ${\rm D}$-module
structure, where
$${\rm D}=\frac{K[\partial_{y_1},\cdots,\partial_{y_n}]}{(\partial_{y_1}^p,\cdots,\partial_{y_n}^p)}
=K[t_1,t_2,\cdots,t_n]= \bigoplus^{n(p-1)}_{\ell=0}{\rm
D}_{\ell}$$ which acts on $R$ through the partial derivations
$\partial_{y_1}$, $\partial_{y_2}$, ..., $\partial_{y_n}$. For any
subspace $V\subset R^{\ell}$, let $\mathbb{L}({\rm
D}_{2\ell-n(p-1)}\cdot V)$ be the linear subspace spanned by ${\rm
D}_{2\ell-n(p-1)}\cdot V\subset R^{n(p-1)-\ell}$. Then we are
reduced to ask if
$${\rm dim}(V)\le {\rm
dim}\,\mathbb{L}({\rm D}_{2\ell-n(p-1)}\cdot V)\quad\text{when\,
$\frac{n(p-1)}{2}\le \ell\le n(p-1)$ ?}$$ Our Lemma \ref{lem3.5}
and Proposition \ref{prop3.6} give an affrmative answer to it.

When $X$ is a smooth projective curve of genus $g\ge 1$, the proof
of theorem is very elementary and simple, which does not need the
more involved arguments of higher dimensional case and shows the
idea of proof best. Thus, although it is a direct corollary of the
general case (Theorem \ref{thm3.7}), we still put its proof in an
independent section. It is also convenient for a reader who is
only interested in the proof for curves.

{\it Acknowledgements:} I would like to thank H\'el\`ene Esnault,
Eckart Viehweg, Hourong Qin, Manfred Lehn, Indranil Biswas, Herbert
Lange, Christian Pauly for their interest and discussions. The proof
of the purely combinatorial Lemma \ref{lem} is due to Fusheng Leng.
I thank him very much for his help. Finally, I would like to thank
the referee for the helpful comments.

\section{The case of curves}
Let $k$ be an algebraically closed field of characteristic $p>0$
and $X$ be a smooth projective curve over $k$. Let $F:X\to X_1$ be
the relative $k$-linear Frobenius morphism, where
$X_1:=X\times_kk$ is the base change of $X/k$ under the Frobenius
$\Spec(k)\to\Spec(k)$. Let $W$ be a vector bundle on $X$ and
$V=F^*(F_*W)$. It is known (\cite[Theorem 5.1]{K}) that $V$ has a
canonical connection $\nabla: V\to V\otimes\Omega^1_X$ with zero
$p$-curvature. In \cite[Section 5]{JRXY}, the authors defined a
canonical filtration \ga{1.1}{0=V_p\subset
V_{p-1}\subset\cdots\subset V_{\ell}\subset
V_{\ell-1}\subset\cdots V_1\subset V_0=V} where $V_1={\rm
ker}(V=F^*F_*W\surj W)$ and \ga{1.2}{V_{\ell+1}={\rm
ker}(V_{\ell}\xrightarrow{\nabla} V\otimes\Omega^1_X\to
V/V_{\ell}\otimes\Omega^1_X).} The following lemma belongs to them
(cf. \cite[Theorem 5.3]{JRXY}).
\begin{lem}\label{lem1.1}
\begin{itemize}\item[(i)] $V_0/V_1\cong W$, $\nabla(V_{\ell+1})\subset V_{\ell}\otimes\Omega^1_X\,$
 for $\ell\ge 1$.
\item[(ii)]
$V_{\ell}/V_{\ell+1}\xrightarrow{\nabla}(V_{\ell-1}/V_{\ell})\otimes\Omega^1_X$
is an isomorphism for $1\le \ell\le p-1$.
\item[(iii)] If $g\ge 2$ and $W$ is semistable, then the canonical
filtration \eqref{1.1} is nothing but the Harder-Narasimhan
filtration.
\end{itemize}
\end{lem}
\begin{proof} (i) follows by the definition, which and (ii) imply
(iii). To prove (ii), let $I_0=F^*F_*\sO_X$, $I_1={\rm
ker}(F^*F_*\sO_X\surj \sO_X)$ and
 \ga{1.3} {I_{\ell+1}={\rm ker}(I_{\ell}\xrightarrow{\nabla}
I_0\otimes\Omega^1_X\surj I_0/I_{\ell}\otimes\Omega^1_X)} which is
the canonical filtration \eqref{1.1} in the case $W=\sO_X$.

(ii) is clearly a local problem, we can assume $X=\Spec(k[[x]])$ and
$W=k[[x]]^{\oplus r}$. Then $V_0:=V=F^*(F_*W)=I_0^{\oplus r}$,
$V_{\ell}=I_{\ell}^{\oplus r}$ and
\ga{1.4}{V_{\ell}/V_{\ell+1}=(I_{\ell}/I_{\ell+1})^{\oplus
r}\xrightarrow{\oplus\nabla}
(I_{\ell-1}/I_{\ell}\otimes\Omega^1_X)^{\oplus
r}=V_{\ell-1}/V_{\ell}\otimes\Omega^1_X.} Thus it is enough to show
that \ga{1.5}{I_{\ell}/I_{\ell+1}\xrightarrow{\nabla}
I_{\ell-1}/I_{\ell}\otimes\Omega^1_X} is an isomorphism. Locally,
$I_0=k[[x]]\otimes_{k[[x^p]]}k[[x]]$ and
\ga{1.6}{\nabla:k[[x]]\otimes_{k[[x^p]]}k[[x]]\to
I_0\otimes_{\sO_X}\Omega^1_X,} where $\nabla(g\otimes f)=g\otimes
f'\otimes{\rm d}x$. The $\sO_X$-module \ga{1.7} {I_1:={\rm
ker}(k[[x]]\otimes_{k[[x^p]]}k[[x]]\surj k[[x]])} has a basis $\{
x^i\otimes 1-1\otimes x^i\}_{1\le i\le p-1}$. Notice that $I_1$ is
also an ideal of the $\sO_X$-algebra
$I_0=k[[x]]\otimes_{k[[x^p]]}k[[x]]$, let $\alpha=x\otimes
1-1\otimes x$, then $\alpha^i\in I_1$. It is easy to see that
$\alpha,\,\alpha^2,\,\ldots,\,\alpha^{p-1}$ is a basis of the
$\sO_X$-module $I_1$ (notice that $\alpha^p=x^p\otimes 1-1\otimes
x^p=0$), and \ga{1.8}
{\nabla(\alpha^{\ell})=-\ell\alpha^{\ell-1}\otimes {\rm d}x.} Thus,
as a free $\sO_X$-module, $I_{\ell}$ has a basis
$\{\alpha^{\ell},\,\alpha^{\ell+1},\ldots,\,\alpha^{p-1}\}$, which
means that $I_{\ell}/I_{\ell+1}$ has a basis $\alpha^{\ell}$,
$(I_{\ell-1}/I_{\ell})\otimes\Omega^1_X$ has a basis
$\alpha^{\ell-1}\otimes{\rm d}x$ and
$\nabla(\alpha^{\ell})=-\ell\alpha^{\ell-1}\otimes{\rm d}x$.
Therefore $\nabla$ induces the isomorphism \eqref{1.5} since
$(\ell,p)=1$, which implies the isomorphism in (ii).
\end{proof}

\begin{thm}\label{thm1.3} Let $X$ be a smooth projective curve of
genus $g\ge 1$. Then $F_*W$ is semi-stable whenever $W$ is
semi-stable. If $g\ge 2$, then $F_*W$ is stable whenever $W$ is
stable.
\end{thm}
\begin{proof} Let $\sE\subset F_*W$ be a nontrivial subbundle and
\ga{1.9}{0\subset V_m\cap F^*\sE\subset\,\cdots\,\subset V_1\cap
F^*\sE\subset V_0\cap F^*\sE=F^*\sE} be the induced filtration.
Let $r_{\ell}={\rm rk}(\frac{V_{\ell}\cap F^*\sE}{V_{\ell+1}\cap
F^*\sE})$ be the ranks of quotients. Then, by the filtration
\eqref{1.9}, we have \ga{1.10}{\mu(F^*\sE)=\frac{1}{{\rm
rk}(F^*\sE)}\sum^m_{\ell=0}r_{\ell}\cdot\mu(\frac{V_{\ell}\cap
F^*\sE}{V_{\ell+1}\cap F^*\sE}).} By Lemma \ref{lem1.1},
$V_{\ell}/V_{\ell+1}\cong W\otimes(\Omega^1_X)^{\otimes \ell}$ is
stable, we have \ga{1.11} {\mu(\frac{V_{\ell}\cap
F^*\sE}{V_{\ell+1}\cap F^*\sE})\le\mu(W)+2(g-1)\ell.} Then, notice
that $\mu(V)=\mu(W)+(p-1)(g-1)$, we have
\ga{1.12}{\mu(F_*W)-\mu(\sE)\ge \frac{2g-2}{p\cdot{\rm
rk}(\sE)}\cdot\sum^m_{\ell=0}(\frac{p-1}{2}-\ell)r_{\ell}} which
becomes equality if and only if the inequalities in \eqref{1.11}
become equalities. It is clear by \eqref{1.12} that
$\mu(F_*W)-\mu(\sE)>0$ if $m\le \frac{p-1}{2}$. Thus we can assume
that $m>\frac{p-1}{2}$, then we can write
 \ga{1.13}
{\sum^m_{\ell=0}(\frac{p-1}{2}-\ell)r_{\ell}
=\sum^{p-1}_{\ell=m+1}(\ell-\frac{p-1}{2})
r_{p-1-\ell}\\+\sum^m_{\ell>\frac{p-1}{2}}(\ell-\frac{p-1}{2})
(r_{p-1-\ell}-r_{\ell})\\\notag\ge
\sum^m_{\ell>\frac{p-1}{2}}(\ell-\frac{p-1}{2})
(r_{p-1-\ell}-r_{\ell}).\notag} On the other hand, since the
isomorphisms
$V_{\ell}/V_{\ell+1}\xrightarrow{\nabla}(V_{\ell-1}/V_{\ell})\otimes\Omega^1_X$
in Lemma \ref{lem1.1} (ii) induce the injections
$$\frac{V_{\ell}\cap F^*\sE}{V_{\ell+1}\cap F^*\sE}\inj
 \frac{V_{\ell-1}\cap F^*\sE}{V_{\ell}\cap
 F^*\sE}\otimes\Omega^1_X$$
we have $r_0\ge r_1\ge \cdots\ge r_{\ell-1}\ge
r_{\ell}\ge\cdots\ge r_m$. Thus
$$\mu(F_*W)-\mu(\sE)\ge
\frac{2g-2}{p\cdot{\rm
rk}(\sE)}\sum\limits^m_{\ell=0}(\frac{p-1}{2}-\ell)r_{\ell}\,\ge
0\,.$$ If $\mu(F_*W)-\mu(\sE)=0$, then \eqref{1.12} and
\eqref{1.13} become equalities. That \eqref{1.12} becomes equality
implies inequalities in \eqref{1.11} become equalities, which
means $r_0=r_1=\cdots=r_m={\rm rk}(W)$. Then that \eqref{1.13}
become equalities implies $m=p-1$. Altogether imply $\sE=F_*W$, we
get contradiction. Hence $F_*W$ is stable whenever $W$ is stable.
\end{proof}

\section{The filtration on higher dimension varieties}

Let $X$ be a smooth projective variety over $k$ of dimension $n$
and $F:X\to X_1$ be the relative $k$-linear Frobenius morphism,
where $X_1:=X\times_kk$ is the base change of $X/k$ under the
Frobenius $\Spec(k)\to\Spec(k)$. Let $W$ be a vector bundle on $X$
and $V=F^*(F_*W)$. We have the straightforward generalization of
the canonical filtration to higher dimensional varieties.

\begin{defn}\label{defn2.1} Let $V_0:=V=F^*(F_*W)$,
$V_1=\ker(F^*(F_*W)\surj W)$ \ga{2.1}
{V_{\ell+1}:=\ker(V_{\ell}\xrightarrow{\nabla} V\otimes_{\sO_X}
\Omega^1_X\to (V/V_{\ell})\otimes_{\sO_X}\Omega^1_X)} where
$\nabla: V\to V\otimes_{\sO_X} \Omega^1_X$ is the canonical
connection (cf. \cite[Theorem 5.1]{K}).
\end{defn}

We first consider the special case $W=\sO_X$ and give some local
descriptions. Let $I_0=F^*(F_*\sO_X)$, $I_1=\ker(F^*F_*\sO_X\surj
\sO_X)$ and \ga{2.2}
{I_{\ell+1}=\ker(I_{\ell}\xrightarrow{\nabla}I_0\otimes_{\sO_X}\Omega^1_X\to
I_0/I_{\ell}\otimes_{\sO_X}\Omega^1_X).}

Locally, let $X=\Spec(A)$, $I_0=A\otimes_{A^p} A$, where $
A=k[[x_1,\cdots,x_n]]$, $A^p=k[[x^p_1,\cdots,x^p_n]]$. Then the
canonical connection $\nabla: I_0\to I_0\otimes\Omega^1_X$ is
locally defined by \ga{2.3} {\nabla(g\otimes_{A^p} f)=\sum_{i=1}^n
(g\otimes_{A^p}\frac{\partial f}{\partial x_i})\otimes_A {\rm
d}x_i}  Notice that $I_0$ has an $A$-algebra structure such that
$I_0=A\otimes_{A^p}A\surj A$ is a homomorphism of $A$-algebras,
its kernel $I_1$ contains elements \ga{2.4}
{\alpha_1^{k_1}\alpha_2^{k_2}\cdots \alpha_n^{k_n},\,\,\,
\text{where $\alpha_i=x_i\otimes_{A^p} 1-1\otimes_{A^p} x_i$,\,\,
$\sum^n_{i=1} k_i\ge 1.$}} Since
$\alpha_i^p=x_i^p\otimes_{A^p}1-1\otimes_{A^p}x_i^p=0$, the set
$\{\alpha_1^{k_1}\cdots \alpha_n^{k_n}\,|\, k_1+\cdots+k_n\ge 1\}$
has $p^n-1$ elements. In fact, we have
\begin{lem}\label{lem2.2} Locally, as
free $A$-modules, we have, for all $\ell\ge 1$, \ga{2.5}
{I_{\ell}=\bigoplus_{k_1+\cdots+k_n\ge
\ell}(\alpha_1^{k_1}\cdots\alpha_n^{k_n})A.}
\end{lem}
\begin{proof} We first prove for $\ell=1$ that
$\{\alpha_1^{k_1}\cdots\alpha_n^{k_n}\,|\,k_1+\cdots+k_n\ge 1\}$ is
a basis of $I_1$ locally. By definition, $I_1$ is locally free of
rank $p^n-1$, thus it is enough to show that as an $A$-module $I_1$
is generated locally by $\{\alpha_1^{k_1}\cdots \alpha_n^{k_n}\,|\,
k_1+\cdots+k_n\ge 1\}$ since it has exactly $p^n-1$ elements.

It is easy to see that as an $A$-module $I_1$ is locally generated
by $$\{x_1^{k_1}\cdots
x_n^{k_n}\otimes_{A^p}1-1\otimes_{A^p}x_1^{k_1}\cdots
x_n^{k_n}\,|\,k_1+\cdots+k_n\ge 1,\,\,0\le k_i\le p-1\,\}.$$ It is
enough to show that any $x_1^{k_1}\cdots
x_n^{k_n}\otimes_{A^p}1-1\otimes_{A^p}x_1^{k_1}\cdots x_n^{k_n}$ is
a linear combination of
$\{\alpha_1^{k_1}\cdots\alpha_n^{k_n}\,|\,k_1+\cdots+k_n\ge 1\}.$
The claim is obvious when $k_1+\cdots+k_n=1$, we consider the case
$k_1+\cdots+k_n>1$. Without loss generality, assume $k_n\ge 1$ and
there are $f_{j_1,\ldots,j_n}\in A$ such that
$$x_1^{k_1}\cdots
x_n^{k_n-1}\otimes_{A^p}1-1\otimes_{A^p}x_1^{k_1}\cdots
x_n^{k_n-1}=\sum_{j_1+\cdots+j_n\ge
1}(\alpha_1^{j_1}\cdots\alpha_n^{j_n})\cdot f_{j_1,\ldots,j_n}.$$
Then we have
$$\aligned &x_1^{k_1}\cdots
x_n^{k_n}\otimes_{A^p}1-1\otimes_{A^p}x_1^{k_1}\cdots
x_n^{k_n}=\sum_{j_1+\cdots+j_n\ge
1}(\alpha_1^{j_1}\cdots\alpha_n^{j_n+1})\cdot
f_{j_1,\ldots,j_n}\\& +\sum_{j_1+\cdots+j_n\ge
1}(\alpha_1^{j_1}\cdots\alpha_n^{j_n})\cdot f_{j_1,\ldots,j_n}x_n
\quad+\quad\alpha_n\cdot (x_1^{k_1}\cdots x_n^{k_n-1}).
\endaligned$$

For $\ell>1$, to prove the lemma, we first show
\ga{2.6}{\nabla(\alpha_1^{k_1}\cdots\alpha_n^{k_n})=-\sum^n_{i=1}k_i
(\alpha_1^{k_1}\cdots\alpha_i^{k_i-1}\cdots\alpha_n^{k_n})
\otimes_A{\rm d}x_i}  Indeed, \eqref{2.6} is true when
$k_1+\cdots+k_n=1$. If $k_1+\cdots+k_n>1$, we assume $k_n\ge 1$
and $\alpha_1^{k_1}\cdots\alpha_n^{k_n-1}=\sum
g_j\otimes_{A^p}f_j$. Then
$$\alpha_1^{k_1}\cdots\alpha_n^{k_n}=\sum_j
x_ng_j\otimes_{A^p}f_j-\sum_jg_j\otimes_{A^p}f_jx_n\,.$$ Use
\eqref{2.3}, straightforward computations show
$$\nabla(\alpha_1^{k_1}\cdots\alpha_n^{k_n})=\alpha_n\nabla(\alpha_1^{k_1}\cdots\alpha_n^{k_n-1})-
(\alpha_1^{k_1}\cdots\alpha_n^{k_n-1})\otimes_A{\rm d}x_n$$ which
implies \eqref{2.6}. Now we can assume the lemma is true for
$I_{\ell-1}$ and recall that
$I_{\ell}=\ker(I_{\ell-1}\xrightarrow{\nabla}I_0\otimes_A\Omega^1_X\surj
(I_0/I_{\ell-1})\otimes_A\Omega^1_X)$. For any
$$\beta=\sum_{k_1+\cdots k_n\ge
\ell-1}(\alpha_1^{k_1}\cdots\alpha_n^{k_n})\cdot
f_{k_1,\ldots,k_n}\in I_{\ell-1},\quad f_{k_1,\ldots,k_n}\in A,$$
by using \eqref{2.6}, we see that $\beta\in I_{\ell}$ if and only
if \ga{2.7}
{\sum_{k_1+\cdots+k_n=\ell-1}(\alpha_1^{k_1}\cdots\alpha_j^{k_j-1}\cdots\alpha_n^{k_n})\cdot
k_jf_{k_1,\ldots,k_n}\,\,\in I_{\ell-1}} for all $1\le j\le n$.
Since $\{\alpha_1^{k_1}\cdots\alpha_n^{k_n}\,|\,k_1+\cdots+k_n\ge
1\}$ is a basis of $I_1$ locally and the lemma is true for
$I_{\ell-1}$, \eqref {2.7} is equivalent to \ga{2.8} {{\rm For} \
{\rm given} \
(k_1,\ldots,k_n) \ {\rm with} \ k_1+\cdots+k_n=\ell-1\\
k_jf_{k_1,\ldots,k_n}=0 \ \ {\rm for} \  {\rm all} \
j=1,\ldots,n\notag} which implies $f_{k_1,\ldots,k_n}=0$ whenever
$k_1+\cdots+k_n=\ell-1$. Thus $I_{\ell}$ is generated by
$\{\alpha_1^{k_1}\cdots\alpha_n^{k_n}\,\,|\, k_1+\cdots+k_n\ge
\ell\,\}$.
\end{proof}

\begin{lem}\label{lem2.3}
\begin{itemize}\item[(i)]$I_{\ell}=0$ when $\ell>n(p-1)$, and $\nabla(I_{\ell+1})\subset I_{\ell}\otimes\Omega^1_X$
for $\ell\ge 1$.
\item[(ii)]
$I_{\ell}/I_{\ell+1}\xrightarrow{\nabla}(I_{\ell-1}/I_{\ell})\otimes\Omega^1_X$
are injective in the category of vector bundles for $1\le \ell\le
n(p-1)$. In particular, their composition \ga{2.9} {\nabla^{\ell}:
I_{\ell}/I_{\ell+1}\to
(I_0/I_1)\otimes_{\sO_X}(\Omega^1_X)^{\otimes
\ell}=(\Omega^1_X)^{\otimes \ell}} is injective in the category of
vector bundles.
\end{itemize}
\end{lem}

\begin{proof} (i) follows from Lemma \ref{2.2} and Definition
\ref{2.1}. (ii) follows from \eqref {2.6}.

\end{proof}

In order to describe the image of $\nabla^{\ell}$ in \eqref{2.9},
we recall a ${\rm GL}(n)$-representation ${\rm T}^{\ell}(V)\subset
V^{\otimes\ell}$ where $V$ is the standard representation of ${\rm
GL}(n)$. Let ${\rm S}_{\ell}$ be the symmetric group of $\ell$
elements with the action on $V^{\otimes\ell}$ by
$(v_1\otimes\cdots\otimes
v_{\ell})\cdot\sigma=v_{\sigma(1)}\otimes\cdots\otimes
v_{\sigma(\ell)}$ for $v_i\in V$ and $\sigma\in{\rm S}_{\ell}$.
Let $e_1,\,\ldots,\,e_n$ be a basis of $V$, for $k_i\ge 0$ with
$k_1+\cdots+k_n=\ell$ define \ga{2.10}
{v(k_1,\ldots,k_n)=\sum_{\sigma\in{\rm S}_{\ell}}(e_1^{\otimes
k_1}\otimes\cdots\otimes e_n^{\otimes k_n})\cdot\sigma }

\begin{defn}\label{defn2.4} Let ${\rm T}^{\ell}(V)\subset V^{\otimes\ell}$ be
the linear subspace generated by all vectors $v(k_1,\ldots,k_n)$
for all $k_i\ge 0$ satisfying $k_1+\cdots+k_n=\ell$. It is clearly
a representation of ${\rm GL}(V)$. If $\sV$ is a vector bundle of
rank $n$,  the subbundle ${\rm T}^{\ell}(\sV)\subset
\sV^{\otimes\ell}$ is defined to be the associated bundle of the
frame bundle of $\sV$ (which is a principal ${\rm GL}(n)$-bundle)
through the representation ${\rm T}^{\ell}(V)$.
\end{defn}

By sending any $e_1^{k_1}e_2^{k_2}\cdots e_n^{k_n}\in {\rm
Sym}^{\ell}(V)$ to $v(k_1,\ldots,k_n)$, we have \ga{2.11}{{\rm
Sym}^{\ell}(V)\surj {\rm T}^{\ell}(V)} which is an isomorphism in
characteristic zero. When ${\rm char}(k)=p>0$, we have
$v(k_1,\ldots,k_n)=0$ if one of $k_1,\,\ldots,\,k_n$ is bigger
than $p-1$. Thus \eqref{2.11} is not injective when $\ell\ge p$,
and ${\rm T}^{\ell}(V)$ is isomorphic to the quotient of ${\rm
Sym}^{\ell}(V)$ by the relations $e_i^p=0$, $1\le i\le n$. In
particular, \ga{2.12} {{\rm T}^{\ell}(V)\cong{\rm Sym}^{\ell}(V)
\quad {\rm when} \quad 0<\ell<p\,} and ${\rm T}^{\ell}(V)=0$ if
$\ell>n(p-1)$. For any $0<\ell\le n(p-1)$, ${\rm T}^{\ell}(V)$ is
a simple representation of highest weight
$$(\overbrace{p-1,\cdots,p-1}^a,\overbrace{b,0,\cdots,0}^{n-a}),\quad
\text{where $\ell=(p-1)a+b, \,0\le b<p-1$}$$ and is called a
`Truncated symmetric power' (cf. \cite{DW}). In next proposition,
we will describe ${\rm T}^{\ell}(V)$ using symmetric powers and
exterior powers. The case of ${\rm GL}(2)$ is extremely simple, it
is a tensor product of symmetric powers and exterior powers. In
general, let $F^*V$ denote the Frobenius twist of the standard
representation $V$ of ${\rm GL}(n)$ through the homomorphism ${\rm
GL}(n)\to {\rm GL}(n)$ ($(a_{ij})_{n\times n}\to
(a_{ij}^p)_{n\times n}$), we have only a resolution of ${\rm
T}^{\ell}(V)$ using symmetric powers of $V$ and exterior powers of
$F^*V$. Fix a basis $e_1$, ... , $e_n$ of $V$, we define the
$k$-linear maps \ga{2.13}{{\rm Sym}^{\ell-q\cdot
p}(V)\otimes_k\bigwedge^q(V)\xrightarrow{\phi} {\rm
Sym}^{\ell-(q-1)\cdot p}(V)\otimes_k\bigwedge^{q-1}(V)} such that
for any $h=f_{\ell-q\cdot p}\otimes e_{k_1}\wedge\cdots\wedge
e_{k_q}$ ($k_1<\cdots<k_q$), we have
\ga{2.14}{\phi(h)=\sum^q_{i=1}(-1)^{i-1}e_{k_i}^pf_{\ell-q\cdot
p}\otimes e_{k_1}\wedge\cdots
\wedge\hat{e}_{k_i}\wedge\cdots\wedge e_{k_q}.}

\begin{prop}\label{prop2.5} (i) When $n=2$, as ${\rm
GL}(2)$-representations, we have
$${\rm T}^{\ell}(V)=\left\{
\begin{array}{llll} {\rm Sym}^{\ell}(V) &\mbox{when $\ell<p$;}\\
{\rm Sym}^{2(p-1)-\ell}(V)\otimes{\rm
det}(V)^{\ell-(p-1)}&\mbox{when $\ell\ge p$}
\end{array}\right.$$

(ii) Let $\ell(p)\ge 0$ be the unique integer such that
$0\le\ell-\ell(p)\cdot p<p$. Then, in the category of ${\rm
GL}(n)$-representations, we have exact sequence
$$\aligned&0\to{\rm Sym}^{\ell-\ell(p)\cdot
p}(V)\otimes_k\bigwedge^{\ell(p)}(F^*V)\xrightarrow{\phi}{\rm
Sym}^{\ell-(\ell(p)-1)\cdot
p}(V)\otimes_k\bigwedge^{\ell(p)-1}(F^*V)\\&\to\cdots\to {\rm
Sym}^{\ell-q\cdot
p}(V)\otimes_k\bigwedge^q(F^*V)\xrightarrow{\phi}{\rm
Sym}^{\ell-(q-1)\cdot
p}(V)\otimes_k\bigwedge^{q-1}(F^*V)\\&\to\cdots\to{\rm
Sym}^{\ell-p}(V)\otimes_k F^*V\xrightarrow{\phi}{\rm
Sym}^{\ell}(V)\to{\rm T}^{\ell}(V)\to 0.
\endaligned$$
\end{prop}

\begin{proof} (i) When $\ell<p$, ${\rm T}^{\ell}(V)={\rm
Sym}^{\ell}(V)$ follows the construction. When $\ell\ge p$, the
simple representation ${\rm T}^{\ell}(V)$ has highest weight
$$(p-1,\ell-p+1)=(2p-2-\ell,0)+(\ell-p+1)\cdot(1,1)$$
where $(2p-2-\ell,0)$ and $(1,1)$ are the highest weights of the
simple representations ${\rm Sym}^{2p-2-\ell}(V)$ and
$\wedge^2(V)={\rm det}(V)$ respectively. Thus
$${\rm T}^{\ell}(V)={\rm Sym}^{2(p-1)-\ell}(V)\otimes{\rm
det}(V)^{\ell-(p-1)}.$$

(ii) The elements $e_1^p,\,e_2^p,\, \ldots,\,e_n^p\in {\rm
Sym}^{\bullet}(V)$ form clearly a regular sequence for ${\rm
Sym}^{\bullet}(V)$, thus the Koszul complex
$K_{\bullet}(e_1^p,\ldots,e_n^p)$ of ${\rm
Sym}^{\bullet}(V)$-modules is a resolution of
$$\frac{{\rm
Sym}^{\bullet}(V)}{(e_1^p,e_2^p,\ldots,e_n^p){\rm
Sym}^{\bullet}(V)}$$ where $K_1={\rm Sym}^{\bullet}(V)\otimes_kV$
with basis $1\otimes_ke_1,\,\ldots,\,1\otimes_ke_n$ and
$K_i=\wedge^iK_1$. Notice $\wedge^iK_1\cong {\rm
Sym}^{\bullet}(V)\otimes_k\wedge^iV$ (as ${\rm
Sym}^{\bullet}(V)$-modules), the sequence in the proposition is
exact in the category of $k$-linear spaces (This was pointed out
by Manfred Lehn).

We only need to show the $k$-linear maps $\phi$ in \eqref{2.13}
are maps of ${\rm GL}(n)$-representations if $\wedge^{\cdot}V$ is
twisted by Frobenius. It is enough to show, for any
$A=(a_{ij})_{n\times n}\in {\rm GL}(n)$ and $h=1\otimes
e_{k_1}\wedge\cdots\wedge e_{k_q}$, that
$$\phi(A\cdot h)=A\cdot\phi(h)$$
To simplify notation, we assume $h=1\otimes e_1\wedge\cdots\wedge
e_q$, then
$$A\cdot h=1\otimes\sum_{k_1<\cdots<k_q}D\left(\begin{array}{cccc}
k_1,k_2,\ldots,k_q\\
1,2,\ldots,q\end{array}\right) e_{k_1}\wedge\cdots\wedge
e_{k_q},\quad \text{where}$$
$$D\left(\begin{array}{cccc}
k_1,k_2,\ldots,k_q\\
1,2,\ldots,q\end{array}\right)=\left|
\begin{array}{cccc}
a^p_{1k_1},a^p_{1k_2},\ldots,a^p_{1k_q}\\
a^p_{2k_1},a^p_{2k_2},\ldots,a^p_{2k_q}\\
\vdots\qquad\vdots\quad\cdots\quad\vdots\\
 a^p_{qk_1},a^p_{qk_2},\ldots,a^p_{qk_q}
\end{array} \right|$$
Then, by definition of $\phi$, we have
$$\aligned&\phi(A\cdot h)=\sum_{k_1<\cdots<k_q}\sum^q_{i=1}(-1)^{i-1}e^p_{k_i}\otimes_k D\left(\begin{array}{cccc}
k_1,k_2,\ldots,k_q\\
1,2,\ldots,q\end{array}\right) e_{k_1}\wedge\cdots\wedge\hat
e_{k_i}\wedge \cdots\wedge e_{k_q}\\&=
\sum_{k_1<\cdots<k_q}\sum^q_{i=1}e^p_{k_i}\otimes_k
\sum^q_{j=1}(-1)^{j-1}a^p_{jk_i}D\left(\begin{array}{cccc}
k_1,\ldots,\hat k_i,\ldots,k_q\\
1,\ldots,\hat j,\ldots,q\end{array}\right)
e_{k_1}\wedge\cdots\wedge\hat e_{k_i}\wedge \cdots\wedge
e_{k_q}\\&=
\sum^q_{j=1}(-1)^{j-1}\sum_{k_1<\cdots<k_q}\sum^q_{i=1}a^p_{jk_i}e^p_{k_i}\otimes_k
D\left(\begin{array}{cccc}
k_1,\ldots,\hat k_i,\ldots,k_q\\
1,\ldots,\hat j,\ldots,q\end{array}\right)
e_{k_1}\wedge\cdots\wedge\hat e_{k_i}\wedge \cdots\wedge
e_{k_q}.\endaligned$$ On the other hand, we will show
$$\aligned&\sum_{k_1<\cdots<k_q}\sum^q_{i=1}a^p_{jk_i}e^p_{k_i}\otimes_k
D\left(\begin{array}{cccc}
k_1,\ldots,\hat k_i,\ldots,k_q\\
1,\ldots,\hat j,\ldots,q\end{array}\right)
e_{k_1}\wedge\cdots\wedge\hat e_{k_i}\wedge \cdots\wedge
e_{k_q}\\&=\left(\sum^n_{i=1}a^p_{ji}e^p_i\right)\otimes_k\left(\sum^n_{i=1}a^p_{1i}e_i\right)\wedge\cdots
\widehat{\left(\sum^n_{i=1}a^p_{ji}e_i\right)}\cdots\wedge\left(\sum^n_{i=1}a^p_{qi}e_i\right)\\&-
\sum_{i_1<\cdots<i_{q-1}}\left(\sum^{q-1}_{k=1}a^p_{ji_k}e^p_{i_k}\right)\otimes_k
D\left(\begin{array}{cccc}
i_1,\ldots,i_{q-1}\\
1,\ldots\hat j\ldots, q\end{array}\right)e_{i_1}\wedge\cdots\wedge
e_{i_{q-1}}
\endaligned$$
and $\sum\limits^q_{j=1}(-1)^{j-1}a^p_{ji_k}\cdot
D\left(\begin{array}{cccc}
i_1,\ldots,i_{q-1}\\
1,\ldots\hat j\ldots, q\end{array}\right)=0$ ($1\le k\le q-1$).
Thus $$\phi(A\cdot h)=A\cdot\phi(h).$$ In fact, the second
equality corresponds to developing a determinant having the
$i_k$-th column repeated. To show the first equality, write
$$\aligned&\left(\sum^n_{i=1}a^p_{ji}e^p_i\right)\otimes_k\left(\sum^n_{i=1}a^p_{1i}e_i\right)\wedge\cdots
\wedge\widehat{\left(\sum^n_{i=1}a^p_{ji}e_i\right)}\wedge\cdots\wedge\left(\sum^n_{i=1}a^p_{qi}e_i\right)\\&=
\sum_{i_1<\cdots<i_{q-1}}\left(\sum^n_{i=1}a^p_{ji}e^p_i
\right)\otimes_k D\left(\begin{array}{cccc}
i_1,\ldots,i_{q-1}\\
1,\ldots\hat j\ldots, q\end{array}\right)e_{i_1}\wedge\cdots\wedge
e_{i_{q-1}}\endaligned.$$ For given $i_1<\cdots<i_{q-1}$, let
$S=\{i_1,\,\ldots,\,i_{q-1}\}$, write
$$\aligned
&\left(\sum^n_{i=1}a^p_{ji}e^p_i \right)\otimes_k
D\left(\begin{array}{cccc}
i_1,\ldots,i_{q-1}\\
1,\ldots\hat j\ldots, q\end{array}\right)e_{i_1}\wedge\cdots\wedge
e_{i_{q-1}}=\\&\sum_{t\notin
S}a^p_{jt}e^p_t\otimes_kD\left(\begin{array}{cccc}
i_1,\ldots,i_{q-1}\\
1,\ldots\hat j\ldots, q\end{array}\right)e_{i_1}\wedge\cdots\wedge
e_{i_{q-1}}+\\&\sum^{q-1}_{k=1}a^p_{ji_k}e^p_{i_k}\otimes_kD\left(\begin{array}{cccc}
i_1,\ldots,i_{q-1}\\
1,\ldots\hat j\ldots, q\end{array}\right)e_{i_1}\wedge\cdots\wedge
e_{i_{q-1}}
\endaligned$$
notice that for any $t\notin S$ there is a unique $k_1<\cdots<k_q$
with $k_i=t$ such that $(k_1,...,\hat
k_i,...,k_q)=(i_1,...,i_{q-1})$, we have
$$\aligned&\sum_{t\notin
S}a^p_{jt}e^p_t\otimes_kD\left(\begin{array}{cccc}
i_1,\ldots,i_{q-1}\\
1,\ldots\hat j\ldots, q\end{array}\right)e_{i_1}\wedge\cdots\wedge
e_{i_{q-1}}=\\&\sum_{k_1<\cdots<k_q}a^p_{jk_i}e^p_{k_i}\otimes_k
D\left(\begin{array}{cccc}
k_1,\ldots,\hat k_i,\ldots,k_q\\
1,\ldots,\hat j,\ldots,q\end{array}\right)
e_{k_1}\wedge\cdots\wedge\hat e_{k_i}\wedge \cdots\wedge
e_{k_q}\endaligned$$ where the summation is taken for all
$k_1<\cdots <k_q$ satisfying $$(k_1,...,\hat
k_i,...,k_q)=(i_1,...,i_{q-1}).$$ Then, taking summation for all
$i_1<\cdots<i_{q-1}$ and exchange the order of two summations, we
got the claimed equality.
\end{proof}

\begin{lem}\label{lem2.6} With the notation in Definition
\ref{defn2.4}, the composition \ga{2.15}
{\nabla^{\ell}:I_{\ell}/I_{\ell+1}\to (\Omega_X^1)^{\otimes\ell}}
of the $\sO_X$-morphisms in Lemma \ref{lem2.3} (ii) has image
${\rm T}^{\ell}(\Omega^1_X)\subset (\Omega_X^1)^{\otimes\ell}$.
\end{lem}
\begin{proof} It is enough to prove the lemma locally. By Lemma
\ref{lem2.2}, $I_{\ell}/I_{\ell+1}$ is locally generated by
\ga{2.16} {\{\alpha_1^{k_1}\cdots\alpha_n^{k_n}\,\,|\,
k_1+\cdots+k_n= \ell\,\}.} By using formula \eqref{2.6} and the
formula of permutations with repeated objects, we have \ga{2.17}
{\nabla^{\ell}(\alpha_1^{k_1}\cdots\alpha_n^{k_n})=(-1)^{\ell}\sum_{\sigma\in{\rm
S}_{\ell}}({\rm d}x_1^{\otimes k_1}\otimes\cdots{\rm
d}x_n^{\otimes k_n})\cdot\sigma } which implies that
$\nabla^{\ell}(I_{\ell}/I_{\ell+1})={\rm T}^{\ell}
(\Omega^1_X)\subset (\Omega_X^1)^{\otimes\ell}$.
\end{proof}

\begin{thm}\label{thm2.7} The filtration defined in Definition
\ref{defn2.1} is \ga{2.18} {0=V_{n(p-1)+1}\subset
V_{n(p-1)}\subset\cdots\subset V_1\subset V_0=V=F^*(F_*W)} which
has the following properties
\begin{itemize}\item[(i)]$\nabla(V_{\ell+1})\subset V_{\ell}\otimes\Omega^1_X$
for $\ell\ge 1$, and $V_0/V_1\cong W$.
\item[(ii)]
$V_{\ell}/V_{\ell+1}\xrightarrow{\nabla}(V_{\ell-1}/V_{\ell})\otimes\Omega^1_X$
are injective morphisms of vector bundles for $1\le \ell\le
n(p-1)$, which induced isomorphisms $$\nabla^{\ell}:
V_{\ell}/V_{\ell+1}\cong W\otimes_{\sO_X}{\rm
T}^{\ell}(\Omega^1_X),\quad 0\le \ell\le n(p-1).$$ The vector
bundle ${\rm T}^{\ell}(\Omega^1_X)$ is suited in the exact
sequence
$$\aligned&0\to{\rm Sym}^{\ell-\ell(p)\cdot p}(\Omega^1_X)\otimes
F^*\Omega_X^{\ell(p)}\xrightarrow{\phi}{\rm
Sym}^{\ell-(\ell(p)-1)\cdot p}(\Omega^1_X)\otimes
F^*\Omega_X^{\ell(p)-1}\\&\to\cdots\to {\rm Sym}^{\ell-q\cdot
p}(\Omega^1_X)\otimes F^*\Omega_X^q\xrightarrow{\phi}{\rm
Sym}^{\ell-(q-1)\cdot p}(\Omega^1_X)\otimes
F^*\Omega_X^{q-1}\\&\to\cdots\to{\rm
Sym}^{\ell-p}(\Omega^1_X)\otimes
F^*\Omega^1_X\xrightarrow{\phi}{\rm Sym}^{\ell}(\Omega^1_X)\to{\rm
T}^{\ell}(\Omega_X^1)\to 0
\endaligned$$
where $\ell(p)\ge 0$ is the integer such that $\ell-\ell(p)\cdot
p<p$.
\end{itemize}
\end{thm}

\begin{proof} It is a local problem to prove the theorem. Thus
$V_{n(p-1)+1}=0$ follows from Lemma \ref{lem2.2}. (i) is nothing
but the definition. (ii) follows from Lemma \ref{lem2.3},
Proposition \ref{prop2.5} and Lemma \ref{lem2.6}.
\end{proof}

\begin{cor}\label{cor2.8} When ${\rm dim}(X)=2$, we have
$$V_{\ell}/V_{\ell+1}=\left\{
\begin{array}{llll} W\otimes{\rm Sym}^{\ell}(\Omega^1_X) &\mbox{when $\ell<p$}\\
W\otimes{\rm
Sym}^{2(p-1)-\ell}(\Omega^1_X)\otimes\omega_X^{\ell-(p-1)}&\mbox{when
$\ell\ge p$}
\end{array}\right.$$
\end{cor}

\begin{proof} It follows from (i) of Proposition \ref{prop2.5}.
\end{proof}

\section{stability in higher dimensional case}

Let $X$ be a smooth projective variety over $k$ of dimension $n$
and ${\rm H}$ a fixed ample divisor on $X$. For a torsion free
sheaf $\sE$ on $X$, we define
$$\mu(\sE)=\frac{c_1(\sE)\cdot{\rm H}^{n-1}}{{\rm rk}(\sE)}.$$

\begin{defn}\label{defn3.1} A torsion free sheaf $\sE$ on $X$ is
called semistable (resp. stable) if, for any $0\neq
\sE'\subset\sE$, we have $$\mu(\sE')\le\mu(\sE)\quad ({\rm
resp.}\,\,\mu(\sE')<\mu(\sE)).$$
\end{defn}

For any torsion free sheaf $E$ on $X$, there is a unique
filtration, the so-called Harder-Narasimhan filtration
$$0=E_0\subset E_1\subset\cdots\subset E_k=E$$
such that $E_i/E_{i-1}$ ($1\le i\le k$) are semistable torsion
free sheaves and
$$\mu_{{\rm max}}(E):=\mu(E_1)>\mu(E_2/E_1)>\cdots>\mu(E_k/E_{k-1}):=\mu_{{\rm
min}}(E).$$ The instability of $E$ was defined as
$${\rm I}(E)=\mu_{{\rm max}}(E)-\mu_{{\rm
min}}(E).$$ Then it is easy to see that for any subsheaf $F\subset
E$ we have \ga{3.1} {\mu(F)-\mu(E)\le {\rm I}(E)=\mu_{{\rm
max}}(E)-\mu_{{\rm min}}(E).}

Let $F:X\to X_1$ be the relative $k$-linear Frobenius morphism and
$W$ a vector bundle of rank $r$ on $X$.

\begin{lem}\label{lem3.2} Let $c_1(\Omega^1_X)=K_X$. Then,
in the chow group ${\rm Ch}(X_1)_{\Q}$,
\ga{3.2}{c_1(F_*W)=\frac{r(p^n-p^{n-1})}{2}K_{X_1}+p^{n-1}c_1(W),\\
\mu(F^*F_*W)=p\cdot\mu(F_*W)= \frac{p-1}{2}K_X\cdot{\rm
H}^{n-1}+\mu(W).\notag}
\end{lem}

\begin{proof} The proof is just an application of Riemann-Roch
theorem. Indeed, by Grothendieck-Riemann-Roch theorem, we have
\ga{3.3} {c_1(F_*W)=\frac{rp^n}{2}K_X+F_*(c_1(W)-\frac{r}{2}K_X).}
We remark here that for any irreducible subvariety $Y\subset X$,
its image $F_X(Y)\subset X$ (under the absolute Frobenius
$F_X:X\to X$) equals to $Y$, and the induced morphism $F_X: Y\to
F_X(Y)=Y$ is nothing but the absolute Frobenius morphism $F_Y:Y\to
Y$ (which has degree $p^{{\rm dim}(Y)}$). In particular,
$F_*(c_1(W)-\frac{r}{2}K_X)=p^{n-1}(c_1(W)-\frac{r}{2}K_{X_1})$
proves \eqref{3.2}. That $\mu(F^*F_*W)=p\cdot\mu(F_*W)$ also
follows from this remark.
\end{proof}

Let $V=F^*F_*W$, recall Theorem \ref{thm2.7}, we have the
canonical filtration \ga{3.4} {0=V_{n(p-1)+1}\subset
V_{n(p-1)}\subset\cdots\subset V_1\subset V_0=V=F^*(F_*W)} with
$V_{\ell}/V_{\ell+1}\cong W\otimes_{\sO_X}{\rm
T}^{\ell}(\Omega^1_X)$.

\begin{lem}\label{lem3.3} With the same notation in Theorem
\ref{thm2.7}, we have \ga{3.5} {c_1({\rm
T}^{\ell}(\Omega^1_X))=\frac{\ell}{n}\left(\sum^{\ell(p)}_{q=0}(-1)^qC^q_n\cdot
C^{\ell-qp}_{n+\ell-q-1}\right)K_X\\
{\rm rk}({\rm
T}^{\ell}(\Omega^1_X))=\sum^{\ell(p)}_{q=0}(-1)^qC^q_n\cdot
C^{\ell-qp}_{n+\ell-q-1}.\notag} In particular, we have $ \mu({\rm
T^{\ell}}(\Omega^1_X))=\frac{\ell}{n}K_X\cdot{\rm H}^{n-1}$.
\end{lem}

\begin{proof} The formula of ${\rm rk}({\rm T}^{\ell}(\Omega^1_X))$
follows directly from the exact sequence in Theorem \ref{thm2.7}
(ii). To compute $c_1({\rm T}^{\ell}(\Omega^1_X))$, we use the
fact that for any vector bundle $E$ of rank $n$, we have \ga{3.6}
{c_1({\rm Sym}^q(E))=C^{q-1}_{n+q-1}\cdot c_1(E)} \ga{3.7} {
c_1(\wedge^qE)=C_{n-1}^{q-1}\cdot c_1(E).} Then, use the exact
sequence in Theorem \ref{thm2.7} (ii) and note that
$$c_1(F^*\Omega^q_X)=p\cdot c_1(\Omega^q_X),$$
we have the formula \eqref{3.5} of $c_1({\rm
T}^{\ell}(\Omega^1_X))$.
\end{proof}

Let $\sE\subset F_*W$ be a nontrivial subsheaf, the canonical
filtration \eqref{3.4} induces the filtration (we assume $V_m\cap
F^*\sE\neq 0$) \ga{3.8}{0\subset V_m\cap
F^*\sE\subset\,\cdots\,\subset V_1\cap F^*\sE\subset V_0\cap
F^*\sE=F^*\sE.}

\begin{lem}\label{lem3.4} In the induced filtration \eqref{3.8}, let
$$\sF_{\ell}:=\frac{V_{\ell}\cap F^*\sE}{V_{\ell+1}\cap
F^*\sE}\subset\frac{V_{\ell}}{V_{\ell+1}}, \qquad r_{\ell}={\rm
rk}(\sF_{\ell}).$$ Then there is an injective morphism
$\sF_{\ell}\xrightarrow{\nabla}\sF_{\ell-1}\otimes\Omega^1_X$ and
\ga{3.9} {\aligned\mu(F_*W)-\mu(\sE)\ge&\frac{K_X\cdot{\rm
H}^{n-1}}{np\cdot{\rm
rk}(\sE)}\sum_{\ell=0}^m(\frac{p-1}{2}n-\ell)\cdot r_{\ell}\\&
-\frac{1}{p}\sum^m_{\ell=0}\frac{r_{\ell}\cdot{\rm I}(W\otimes{\rm
T}^{\ell}(\Omega^1_X))}{{\rm rk}(\sE)}\endaligned} the equality
holds if and only if equalities hold in the inequalities \ga{3.10}
{\mu(\sF_{\ell})-\mu(V_{\ell}/V_{\ell+1})\le {\rm I}(W\otimes{\rm
T}^{\ell}(\Omega^1_X))\quad (0\le \ell\le m).}
\end{lem}

\begin{proof} The injective morphisms $V_{\ell}/V_{\ell+1}\xrightarrow{\nabla}(V_{\ell-1}/V_{\ell})\otimes\Omega^1_X$
in Theorem \ref{thm2.7} (ii) induces clearly the injective
morphisms
$$\sF_{\ell}\xrightarrow{\nabla}\sF_{\ell-1}\otimes\Omega^1_X,\qquad \ell=1,\,\ldots,\,m.$$

To show \eqref{3.9}, note
$\mu(F_*W)-\mu(\sE)=\frac{1}{p}(\mu(F^*F_*W)-\mu(F^*\sE))$ and
$$\mu(F^*\sE)=\frac{1}{{\rm
rk}(\sE)}\sum_{\ell=0}^mr_{\ell}\cdot\mu(\sF_{\ell}),$$ using
Lemma \ref{3.2},
we have \ga{3.11} {\mu(F^*F_*W)-\mu(F^*\sE)=\\
\frac{1}{{\rm
rk}(\sE)}\sum^m_{\ell=0}r_{\ell}\left(\frac{p-1}{2}K_X\cdot{\rm
H}^{n-1}+\mu(W)-\mu(\sF_{\ell})\right).\notag} For
$\sF_{\ell}\subset V_{\ell}/V_{\ell+1}=W\otimes{\rm
T}^{\ell}(\Omega^1_X)$ ($0\le \ell\le m$), using Lemma \ref{3.3},
\ga{3.12} {\mu(\sF_{\ell})\le\mu(W)+\frac{\ell}{n}K_X\cdot{\rm
H}^{n-1}+{\rm I}(W\otimes{\rm T}^{\ell}(\Omega^1_X)).} Substitute
\eqref{3.12} into \eqref{3.11}, one get \eqref{3.9} and the
equality holds if and only if all of inequalities \eqref{3.12}
become equalities.

\end{proof}

Let $K$ be a field of characteristic $p>0$, consider the
$K$-algebra
$$R=\frac{K[y_1,\cdots,y_n]}{(y_1^p,\cdots,y_n^p)}=\bigoplus^{n(p-1)}_{\ell=0}R^{\ell},$$
where $R^{\ell}$ is the $K$-linear space generated by
$$\{\,y_1^{k_1}\cdots y_n^{k_n}\,|\,k_1+\cdots+k_n=\ell,\quad 0\le k_i\le
p-1\,\}.$$ The polynomial ring ${\rm P}=
K[\partial_{y_1},\cdots,\partial_{y_n}]$ acts on $R$ through
partial derivations, which induces a ${\rm P}$-module structure on
$R$. Note that $\partial^p_{y_i}$ ($i=1,2,\ldots,n$) act on $R$
trivially, the ${\rm P}$-module structure is in fact a ${\rm
D}$-module, where
$${\rm D}=\frac{K[\partial_{y_1},\cdots,\partial_{y_n}]}{(\partial_{y_1}^p,\cdots,\partial_{y_n}^p)}
=K[t_1,t_2,\cdots,t_n]= \bigoplus^{n(p-1)}_{\ell=0}{\rm
D}_{\ell},$$ where ${\rm D}_{\ell}$ is the linear space of degree
$\ell$ homogeneous elements and $t_1, \, t_2,\,\ldots\,, t_n$ are
the classes of $\partial_{y_1},
\,\partial_{y_2},\,\ldots\,,\partial_{y_n}$.

\begin{lem}\label{lem3.5} Let $V\subset {\rm
D}_{\ell}$ be a linear subspace. Then, when
$\ell\le\frac{n(p-1)}{2}$, there is a basis $\{d_i\in V\}$ of $V$
and monomials $\{\delta_i \in{\rm D}_{n(p-1)-2\ell}\}$ such that
$\{\delta_id_i\in {\rm D}_{n(p-1)-\ell}\}$ are linearly independent.
\end{lem}

\begin{proof} We reduce firstly the lemma to the case
when $V$ has a basis of monomials. Define the Lexicographic order on
the set of monomials of ${\rm D}_{\ell}$, ${\rm D}_{n(p-1)-\ell}$
respectively. For any $v\in {\rm D}_{\ell}$, one can write uniquely
$$v=\lambda_vm_v+\sum_{m>m_v}\lambda_mm$$
where $0\neq \lambda_v,\, \lambda_m\in K$, $m_v$ and $m$ are
monomials of ${\rm D}_{\ell}$.

Let ${\rm dim}(V)=s$, then it is easy to see that there is a basis
$$d_i=\lambda_im_i+\sum_{m>m_i}\lambda_{i,m}m\,,\,\,\lambda_i\neq 0,\quad (1\le i\le
s)$$ of $V$ such that $\{m_1,\ldots, m_s\}$ are different monomials
of ${\rm D}_{\ell}$. If there are monomials $\{\delta_i\in {\rm
D}_{n(p-1)-2\ell}\}_{1\le i\le s}$ such that $\{\delta_im_i\in{\rm
D}_{n(p-1)-\ell}\}_{1\le i\le s}$ are different monomials, then we
claim that
$$\{\delta_id_i\in{\rm D}_{n(p-1)-\ell}\}_{1\le i\le s}$$
are linearly independent. To prove the claim, we only remark that
for any monomials $m, \,m'\in {\rm D}_{\ell}$ and monomial
$\delta\in {\rm D}_{n(p-1)-2\ell}$, we have
$$m<m' \Rightarrow \delta m<\delta m'\quad\text{whenever $\delta
m$, $\delta m'$ are nonzero}.$$ Thus we have
$$\delta_id_i=\lambda_i\delta_im_i+\sum_{\delta_im>\delta
m_i}\lambda_{i,m}\delta_im\quad (1\le i\le s),$$ which are linearly
independent.

If we identify the set of monomials of ${\rm D}_{\ell}$ with the set
$$M^{\ell}=\{\,v=(v_1,\ldots,v_n)\,|0\le v_i\le
p-1\,\,(1\le i\le n),\,\,\,\sum^n_{i=1}v_i=\ell\,\}.$$ Then the
lemma is equivalent to the existence of an injective map
$$\varphi:M^{\ell}\to M^{n(p-1)-\ell}$$
such that for any $v\in M^{\ell}$, we have $v\le\varphi(v)$:
$v_i\le \varphi(v)_i$ \, ($1\le i\le n$). The existence of
$\varphi$ is a special case of the following lemma.
\end{proof}

For any $(a_1,\ldots,a_n)\in \Z^n_{\ge 0}$, let
$M^{\ell}_n(a_1,\ldots,a_n)$ be the set
$$\{\,v=(v_1,\ldots,v_n)\,|0\le
v_i\le a_i\,\,(1\le i\le n),\,\,\,\sum^n_{i=1}v_i=\ell\,\}.$$ For
any $v\in M^{\ell}_n(a_1,\ldots,a_n)$ and $v'\in
M^{\ell\,'}_n(a_1,\ldots,a_n)$, by $v\le v'$, we mean $v_i\le v_i'$
($1\le i\le n$). Then we have the following lemma, its proof was
suggested by Fusheng Leng

\begin{lem}\label{lem} Let $\sigma=\sum\limits^n_{i=1}a_i$. Then, when $\ell\le\frac{1}{2}\sigma$,
there exists an injective map
$\varphi:M^{\ell}_n(a_1,\ldots,a_n)\to
M^{\sigma-\ell}_n(a_1,\ldots,a_n)$ such that
$$v\le\varphi(v)\,\,,\quad \forall \,\,v\in
M^{\ell}_n(a_1,\ldots,a_n).$$
\end{lem}

\begin{proof} The strategy of proof is to do induction for $n$ and
$\sigma$. The lemma is clearly true when $n=1$. Assume the lemma
is true for $n-1$. To show the lemma for $n$, we do induction for
$\sigma$. The lemma is trivially true for any $n$ when $\sigma=1$.
Thus we can assume $n\ge 2$ and $\sigma\ge 2$.

Without loss of generality, we assume $a_{n-1}>0$ and $a_n>0$. Let
$$S^{\ell}=\{\,v\in
M^{\ell}_n(a_1,\ldots,a_n)\,|\,v_{n-1}=a_{n-1}\,\,{\rm
or}\,\,v_n=0\,\},$$
$$S^{\sigma-\ell}=\{\,v\in
M^{\sigma-\ell}_n(a_1,\ldots,a_n)\,|\,v_{n-1}=a_{n-1}\,\,{\rm
or}\,\,v_n=0\,\},$$ $C^{\ell}=M^{\ell}_n(a_1,\ldots,a_n)\setminus
S^{\ell}$ and
$C^{\sigma-\ell}=M^{\sigma-\ell}_n(a_1,\ldots,a_n)\setminus
S^{\sigma-\ell}$. We will show the existence of injective maps
$$\varphi_1:S^{\ell}\to S^{\sigma-\ell}\,,\quad
\varphi_2:C^{\ell}\to C^{\sigma-\ell}$$ with $v\le\varphi_1(v)$,
$v\le\varphi_2(v)$ ($\forall\,\,v\in S^{\ell}$, $\forall\,\,v\in
C^{\ell}$) by induction of $n$, $\sigma$ respectively. In order to
use the induction, we identify $S^{\ell}$ (resp.
$S^{\sigma-\ell}$) with $M^{\ell}_{n-1}(a_1,\ldots,a_{n-1}+a_n)$
(resp. $M^{\sigma-\ell}_{n-1}(a_1,\ldots,a_{n-1}+a_n)$) by
$$f_{\ell}: S^{\ell}\to M^{\ell}_{n-1}(a_1,\ldots,a_{n-1}+a_n),\quad
f_{\ell}(v)=(v_1,\ldots,v_{n-2},v_{n-1}+v_n)$$ (resp.
$f_{\sigma-\ell}:S^{\sigma-\ell}\to
M^{\sigma-\ell}_{n-1}(a_1,\ldots,a_{n-1}+a_n)$). Indeed,
$f_{\ell}$ (resp. $f_{\sigma-\ell}$) is a bijective map. To see
the injectivity of $f_{\ell}$, if $f_{\ell}(v)=f_{\ell}(v')$, then
$v_i=v_i'$ ($1\le i\le n-2$) and $v_{n-1}+v_n=v'_{n-1}+v'_n$. We
claim that $v_{n-1}+v_n=v'_{n-1}+v'_n$ implies $v_n=v'_n$ (thus
$v_{n-1}=v'_{n-1}$) since $v,\,v'\in S^{\ell}$. Indeed, if $v_n=0$
then $v'_n=0$, otherwise $v'_{n-1}=a_{n-1}$ (by definition of
$S^{\ell}$) and $v_{n-1}=a_{n-1}+v'_n>a_{n-1}$ (a contradiction to
the definition of $M^{\ell}_n(a_1,\ldots,a_n)$). Similarly,
$v'_n=0$ implies $v_n=0$. If both $v_n$ and $v'_n$ are not zero,
by definition of $S^{\ell}$, $v_{n-1}=a_{n-1}=v'_{n-1}$, thus
$v_n=v'_n$. To see it being surjective, for any $w\in
M^{\ell}_{n-1}(a_1,\ldots,a_{n-1}+a_n)$, notice that $w_i\le a_i$
($1\le i\le n-2$) and $w_{n-1}\le a_{n-1}+a_n$, we define
$$
v=\begin{cases} (w_1,\ldots,w_{n-2},w_{n-1},0) & \text{if
$w_{n-1}\le a_{n-1}$}\\ (w_1,\ldots,w_{n-2}, a_{n-1},
w_{n-1}-a_{n-1}) & \text{if $w_{n-1}>a_{n-1}$}\end{cases}
$$
then $v\in S^{\ell}$ such that $f_{\ell}(v)=w$. Similarly,
$f_{\sigma-\ell}$ is bijective.

By the inductive assumption for $n$, there exists an injective map
$$\psi_1:M^{\ell}_{n-1}(a_1,\ldots,a_{n-1}+a_n)\to
M^{\sigma-\ell}_{n-1}(a_1,\ldots,a_{n-1}+a_n)$$ such that
$v\le\psi_1(v)$ ($\forall\,\,v\in
M^{\ell}_{n-1}(a_1,\ldots,a_{n-1}+a_n)$). Then, we define
$$\varphi_1=f^{-1}_{\sigma-\ell}\cdot\psi_1\cdot
f_{\ell}:\,\,S^{\ell}\to S^{\sigma-\ell}.$$ For any
$v=(v_1,\ldots,v_n)\in S^{\ell}$, we need to show
$v\le\varphi_1(v)$. Let
$$\psi_1(f_{\ell}(v))=(w_1,\ldots,w_{n-2},w_{n-1})
\in M^{\sigma-\ell}_{n-1}(a_1,\ldots,a_{n-1}+a_n).$$ Then $v_i\le
w_i$ ($1\le i\le n-2$), $v_{n-1}+v_n\le w_{n-1}$ and
$$
\varphi_1(v)=\begin{cases} (w_1,\ldots,w_{n-2},w_{n-1},0) &
\text{if $w_{n-1}\le a_{n-1}$}\\ (w_1,\ldots,w_{n-2}, a_{n-1},
w_{n-1}-a_{n-1}) & \text{if $w_{n-1}>a_{n-1}$}\end{cases}
$$ by the definition of $f_{\ell}$, $\psi_1$ and
$f_{\sigma-\ell}$. Thus $v_i\le \varphi_1(v)_i$ ($1\le i\le n-2$).
We still need to check $v_{n-1}\le \varphi_1(v)_{n-1}$ and
$v_n\le\varphi_1(v)_n$. If $v_n=0$ (thus $v_n\le\varphi_1(v)_n$),
then $v_{n-1}\le w_{n-1}$ (since $v_{n-1}+v_n\le w_{n-1}$), thus
$$v_{n-1}\le {\rm min}\{w_{n-1},a_{n-1}\}\le \varphi_1(v)_{n-1}.$$
If $v_n\neq 0$, by the definition of $S^{\ell}$,
$v_{n-1}=a_{n-1}$, which implies
$$a_{n-1}<a_{n-1}+v_n=v_{n-1}+v_n\le w_{n-1}.$$
Thus $\varphi_1(v)_{n-1}=a_{n-1}$ and
$\varphi_1(v)_n=w_{n-1}-a_{n-1}=w_{n-1}-v_{n-1}\ge v_n.$

Next we construct the injective map $\varphi_2: C^{\ell}\to
C^{\sigma-\ell}$ by using induction for $\sigma$. By the
definition of $C^{\ell}$ and $C^{\sigma-\ell}$, we have
$$C^{\ell}=\{v\in
M^{\ell}_n(a_1,\ldots,a_n)\,|\,v_{n-1}\le
a_{n-1}-1,\,\,v_n\ge1\,\}$$
$$C^{\sigma-\ell}=\{v'\in
M^{\sigma-\ell}_n(a_1,\ldots,a_n)\,|\,v'_{n-1}\le
a_{n-1}-1,\,\,v'_n\ge1\,\}.$$ Let
$\bar\sigma=a_1+\cdots+a_{n-2}+(a_{n-1}-1)+(a_n-1)=\sigma-2$ and
$\bar\ell=\ell-1$, we have the following clear identifications
$$\pi_{\ell}:\,C^{\ell}\to
M^{\bar\ell}_n(a_1,\ldots,a_{n-2},a_{n-1}-1,a_n-1)$$
$$\pi_{\sigma-\ell}:\,C^{\sigma-\ell}\to
M^{\bar\sigma-\bar\ell}_n(a_1,\ldots,a_{n-2},a_{n-1}-1,a_n-1)$$
where $\pi_{\ell}(v)=(v_1,\ldots,v_{n-1}, v_n-1)$,
$\pi_{\sigma-\ell}(v')=(v'_1,\ldots,v'_{n-1}, v'_n-1)$. Notice
that $\bar\ell\le\frac{1}{2}\bar\sigma$, by induction for
$\sigma$, there exists an injective map
$$\psi_2:M^{\bar\ell}_n(a_1,\ldots,a_{n-2},a_{n-1}-1,a_n-1)\to
M^{\bar\sigma-\bar\ell}_n(a_1,\ldots,a_{n-2},a_{n-1}-1,a_n-1)$$
such that $v\le\psi_2(v)$ for any $v\in
M^{\bar\ell}_n(a_1,\ldots,a_{n-2},a_{n-1}-1,a_n-1)$. Let
$$\varphi_2=\pi_{\sigma-\ell}^{-1}\cdot\psi_2\cdot\pi_{\ell}:\,\,C^{\ell}\to
C^{\sigma-\ell}.$$ For any $v=(v_1,\ldots,v_n)\in C^{\ell}$, we
have to check that $v\le\varphi_2(v)$. Let
$$\psi_2(\pi_{\ell}(v))=(w_1,\ldots,w_n)\in
M^{\bar\sigma-\bar\ell}_n(a_1,\ldots,a_{n-2},a_{n-1}-1,a_n-1),$$
then $v_i\le w_i$ ($1\le i\le n-1$), $v_n-1\le w_n$ and
$$\varphi_2(v)=(w_1,\ldots,w_{n-1},w_n+1)\in C^{\sigma-\ell}$$
by the definition of $\pi_{\ell}$, $\psi_2$ and
$\pi_{\sigma-\ell}$. Thus $v_n\le w_n+1=\varphi_2(v)_n$ and we
have shown the lemma.

\end{proof}

\begin{prop}\label{prop3.6} Let $V\subset R^{\ell}$ be a linear
subspace,
$\mathbb{L}({\rm D}_{2\ell-n(p-1)}\cdot V)$ be the linear subspace
generated by ${\rm D}_{2\ell-n(p-1)}\cdot V\subset
R^{n(p-1)-\ell}$. Then,
$${\rm dim}(V)\le {\rm
dim}\,\mathbb{L}({\rm D}_{2\ell-n(p-1)}\cdot V)\quad\text{when\,
$\frac{n(p-1)}{2}\le \ell\le n(p-1)$.}$$
\end{prop}

\begin{proof} Let $\omega=y_1^{p-1}y_2^{p-1}\cdots y_n^{p-1}\in
R^{n(p-1)}$. Then the ${\rm D}$-module structure on $R$ induces
surjective morphisms \ga{3.13}{\phi_{\ell}:{\rm
D}_{\ell}\xrightarrow{\cdot\omega} R^{n(p-1)-\ell}} of linear
spaces for any $0\le \ell\le n(p-1)$. They must be isomorphisms
since ${\rm dim}({\rm D}_{\ell})={\rm dim}(R^{n(p-1)-\ell})$. To
show the equality of dimensions, it is enough to show
$${\rm dim}({\rm D}_{\ell})\ge{\rm dim}(R^{n(p-1)-\ell})={\rm
dim}({\rm D}_{n(p-1)-\ell})\ge{\rm dim}(R^{\ell})={\rm dim}({\rm
D}_{\ell}).$$ The two inequalities hold because we have the
surjective homomorphisms $\phi_{\ell}$ and $\phi_{n(p-1)-\ell}$.
The two equalities hold because
$$\bigoplus^{n(p-1)}_{\ell=0}R^{\ell}=R\cong{\rm D}=\bigoplus^{n(p-1)}_{\ell=0}{\rm D}_{\ell}$$
as (graded) $K$-algebras. In particular, \ga{3.14}
{\phi_{n(p-1)-\ell}:{\rm D}_{n(p-1)-\ell}\to R^{\ell}\,,\quad
\phi_{\ell}:{\rm D}_{\ell}\to R^{n(p-1)-\ell}} are isomorphisms.
Since $0\le \bar{\ell}=n(p-1)-\ell\le \frac{n(p-1)}{2}$, we can use
Lemma \ref{3.5} for $V'=\phi_{n(p-1)-\ell}^{-1}(V)\subset {\rm
D}_{\bar{\ell}}={\rm D}_{n(p-1)-\ell}\,\,$ , thus there is a basis
$\{d_i\in V'\}_{1\le i\le s}$ and monomials $\{\delta_i\in {\rm
D}_{n(p-1)-2\bar\ell}={\rm D}_{2\ell-n(p-1)}\}_{1\le i\le s}$ such
that $\{\delta_id_i\in {\rm D}_{n(p-1)-\bar\ell}={\rm
D}_{\ell}\}_{1\le i\le s}$ are linearly independent. Thus
$$\{\phi_{\ell}(\delta_id_i)=\delta_i(d_i\omega)\in {\rm D}_{2\ell-n(p-1)}\cdot V\subset R^{n(p-1)-\ell}\}_{1\le i\le s}$$
are linearly independent, where $s={\rm dim}(V')={\rm dim}(V)$. We
have proven the proposition.
\end{proof}

Let $X$ be an irreducible smooth projective variety of dimension
$n$ over an algebraically closed field $k$ with ${\rm
char}(k)=p>0$. For any vector bundle $W$ on $X$, let
$${\rm I}(W,X)={\rm max}\{{\rm I}(W\otimes{\rm
T}^{\ell}(\Omega^1_X))\,|\,\, 0\le \ell\le n(p-1)\,\}$$ be the
maximal value of instabilities ${\rm I}(W\otimes{\rm
T}^{\ell}(\Omega^1_X))$.

\begin{thm}\label{thm3.7} When
$K_X\cdot{\rm H}^{n-1}\ge 0$, we have, for any $\sE\subset F_*W$,
\ga{3.15}{\mu(F_*W)-\mu(\sE)\ge -\frac{{\rm I}(W,X)}{p}.} In
particular, if $W\otimes{\rm T}^{\ell}(\Omega^1_X)$, $0\le \ell\le
n(p-1)$, are semistable, then $F_*W$ is semistable. Moreover, if
$K_X\cdot{\rm H}^{n-1}>0$, the stability of the bundles
$W\otimes{\rm T}^{\ell}(\Omega^1_X)$, $0\le \ell\le n(p-1)$,
implies the stability of $F_*W$.
\end{thm}

\begin{proof} Since $K_X\cdot{\rm H}^{n-1}\ge 0$, by the inequality \eqref{3.9} in
Lemma \ref{3.4} (see also the notation in \eqref{3.8} and the
lemma), it is enough to show
$$\sum_{\ell=0}^m(\frac{n(p-1)}{2}-\ell)r_{\ell}\ge 0.$$
If $m\le \frac{n(p-1)}{2}$, it is clear. If $m> \frac{n(p-1)}{2}$,
then we have \ga{3.16}
{\sum_{\ell=0}^m(\frac{n(p-1)}{2}-\ell)r_{\ell}=\sum^{n(p-1)}_{\ell=m+1}
(\ell-\frac{n(p-1)}{2})
r_{n(p-1)-\ell}\\+\sum^m_{\ell\,>\frac{n(p-1)}{2}}(\ell-\frac{n(p-1)}{2})
(r_{n(p-1)-\ell}-r_{\ell}).\notag}

We will use Proposition \ref{prop3.6} to show that
$$r_{\ell}\le r_{n(p-1)-\ell}\quad\text{when\,
$\frac{n(p-1)}{2}\le\ell\le n(p-1)$.}$$ It is clearly a local
problem, we can consider all of the torsion free sheaves as vector
spaces over the function field $K=k(X)$ of $X$. Without loss of
generality, we assume ${\rm rk}(W)=1$. Then, from the discussions
in Section 3, we know that $V_{\ell}/V_{\ell+1}\cong {\rm
T}^{\ell}(\Omega^1_X)$ ($0\le \ell\le n(p-1)$) are precisely
isomorphic to $R^{\ell}$ ($0\le \ell\le n(p-1)$) in Proposition
\ref{prop3.6}. Since the morphisms
$V_{\ell}/V_{\ell+1}\xrightarrow{\nabla}(V_{\ell-1}/V_{\ell})\otimes\Omega^1_X$
induce morphisms
$\sF_{\ell}\xrightarrow{\nabla}\sF_{\ell-1}\otimes\Omega^1_X$, by
the formula \eqref{2.6}, we have
$${\rm D}_{2\ell-n(p-1)}\cdot
\sF_{\ell}\subset\sF_{n(p-1)-\ell}\,.$$ Then, by Proposition
\ref{prop3.6}, $r_{\ell}={\rm dim}(\sF_{\ell})\le{\rm
dim}\,\mathbb{L}({\rm D}_{2\ell-n(p-1)}\cdot \sF_{\ell})\,$, we
have $r_{\ell}\le r_{n(p-1)-\ell}$, thus \eqref{3.15}.

If the bundles $W\otimes{\rm T}^{\ell}(\Omega^1_X)$ ($0\le \ell\le
n(p-1)$) are stable, then
$$\mu(F_*W)-\mu(\sE)\ge 0.$$
It becomes equality if and only if inequalities \eqref{3.10}
become equalities and
$\sum\limits^m_{\ell=0}(\frac{n(p-1)}{2}-\ell)r_{\ell}=0$. Thus
$m>\frac{n(p-1)}{2}$ and each term in \eqref{3.16} must be zero
(since $K_X\cdot{\rm H}^{n-1}>0$), which forces $m=n(p-1)$. Then
the fact that inequalities \eqref{3.10} become equalities implies
$\sE=F_*W$.
\end{proof}

\begin{cor}\label{cor3.8} Let $X$ be a smooth projective variety of ${\rm dim}(X)=n$, whose
canonical divisor $K_X$ satisfies $K_X\cdot {\rm H}^{n-1}\ge 0$.
Then $${\rm I}(F_*W)\le p^{n-1}{\rm rk}(W)\,{\rm I}(W,X).$$
\end{cor}

\begin{proof} It is just Theorem \ref{thm3.7} plus the following
trivial remark: For any vector bundle $E$, if there is a constant
$\lambda$ satisfying $\mu(E')-\mu(E)\le \lambda$ for any
$E'\subset E$. Then ${\rm I}(E)\le {\rm rk}(E)\lambda$.
\end{proof}

\bibliographystyle{plain}

\renewcommand\refname{References}

\end{document}